\documentclass{article}

\usepackage{epsfig}
\usepackage{latexsym, amsfonts, amssymb}
\newdimen\proofrulebreadth \proofrulebreadth=.05em
\newdimen\proofdotseparation \proofdotseparation=1.25ex
\newdimen\proofrulebaseline \proofrulebaseline=2ex
\newcount\proofdotnumber \proofdotnumber=3
\let\then\relax
\def\hfi{\hskip0pt plus.0001fil}
\mathchardef\squigto="3A3B
\newif\ifinsideprooftree\insideprooftreefalse
\newif\ifonleftofproofrule\onleftofproofrulefalse
\newif\ifproofdots\proofdotsfalse
\newif\ifdoubleproof\doubleprooffalse
\let\wereinproofbit\relax
\newdimen\shortenproofleft
\newdimen\shortenproofright
\newdimen\proofbelowshift
\newbox\proofabove
\newbox\proofbelow
\newbox\proofrulename
\def\shiftproofbelow{\let\next\relax\afterassignment\setshiftproofbelow\dimen0 }
\def\shiftproofbelowneg{\def\next{\multiply\dimen0 by-1 }%
\afterassignment\setshiftproofbelow\dimen0 }
\def\setshiftproofbelow{\next\proofbelowshift=\dimen0 }
\def\setproofrulebreadth{\proofrulebreadth}

\def\prooftree{
\ifnum  \lastpenalty=1
\then   \unpenalty
\else   \onleftofproofrulefalse
\fi
\ifonleftofproofrule
\else   \ifinsideprooftree
        \then   \hskip.5em plus1fil
        \fi
\fi
\bgroup
\setbox\proofbelow=\hbox{}\setbox\proofrulename=\hbox{}%
\let\justifies\proofover\let\leadsto\proofoverdots\let\Justifies\proofoverdbl
\let\using\proofusing\let\[\prooftree
\ifinsideprooftree\let\]\endprooftree\fi
\proofdotsfalse\doubleprooffalse
\let\thickness\setproofrulebreadth
\let\shiftright\shiftproofbelow \let\shift\shiftproofbelow
\let\shiftleft\shiftproofbelowneg
\let\ifwasinsideprooftree\ifinsideprooftree
\insideprooftreetrue
\setbox\proofabove=\hbox\bgroup$\displaystyle 
\let\wereinproofbit\prooftree
\shortenproofleft=0pt \shortenproofright=0pt \proofbelowshift=0pt
\onleftofproofruletrue\penalty1
}

\def\eproofbit{
\ifx    \wereinproofbit\prooftree
\then   \ifcase \lastpenalty
        \then   \shortenproofright=0pt  
        \or     \unpenalty\hfil         
        \or     \unpenalty\unskip       
        \else   \shortenproofright=0pt  
        \fi
\fi
\global\dimen0=\shortenproofleft
\global\dimen1=\shortenproofright
\global\dimen2=\proofrulebreadth
\global\dimen3=\proofbelowshift
\global\dimen4=\proofdotseparation
\global\count255=\proofdotnumber
$\egroup  
\shortenproofleft=\dimen0
\shortenproofright=\dimen1
\proofrulebreadth=\dimen2
\proofbelowshift=\dimen3
\proofdotseparation=\dimen4
\proofdotnumber=\count255
}

\def\proofover{
\eproofbit 
\setbox\proofbelow=\hbox\bgroup 
\let\wereinproofbit\proofover
$\displaystyle
}%
\def\proofoverdbl{
\eproofbit 
\doubleprooftrue
\setbox\proofbelow=\hbox\bgroup 
\let\wereinproofbit\proofoverdbl
$\displaystyle
}%
\def\proofoverdots{
\eproofbit 
\proofdotstrue
\setbox\proofbelow=\hbox\bgroup 
\let\wereinproofbit\proofoverdots
$\displaystyle
}%
\def\proofusing{
\eproofbit 
\setbox\proofrulename=\hbox\bgroup 
\let\wereinproofbit\proofusing
\kern0.3em$
}

\def\endprooftree{
\eproofbit 
  \dimen5 =0pt
\dimen0=\wd\proofabove \advance\dimen0-\shortenproofleft
\advance\dimen0-\shortenproofright
\dimen1=.5\dimen0 \advance\dimen1-.5\wd\proofbelow
\dimen4=\dimen1
\advance\dimen1\proofbelowshift \advance\dimen4-\proofbelowshift
\ifdim  \dimen1<0pt
\then   \advance\shortenproofleft\dimen1
        \advance\dimen0-\dimen1
        \dimen1=0pt
        \ifdim  \shortenproofleft<0pt
        \then   \setbox\proofabove=\hbox{%
                        \kern-\shortenproofleft\unhbox\proofabove}%
                \shortenproofleft=0pt
        \fi
\fi
\ifdim  \dimen4<0pt
\then   \advance\shortenproofright\dimen4
        \advance\dimen0-\dimen4
        \dimen4=0pt
\fi
\ifdim  \shortenproofright<\wd\proofrulename
\then   \shortenproofright=\wd\proofrulename
\fi
\dimen2=\shortenproofleft \advance\dimen2 by\dimen1
\dimen3=\shortenproofright\advance\dimen3 by\dimen4
\ifproofdots
\then
        \dimen6=\shortenproofleft \advance\dimen6 .5\dimen0
        \setbox1=\vbox to\proofdotseparation{\vss\hbox{$\cdot$}\vss}%
        \setbox0=\hbox{%
                \advance\dimen6-.5\wd1
                \kern\dimen6
                $\vcenter to\proofdotnumber\proofdotseparation
                        {\leaders\box1\vfill}$%
                \unhbox\proofrulename}%
\else   \dimen6=\fontdimen22\the\textfont2 
        \dimen7=\dimen6
        \advance\dimen6by.5\proofrulebreadth
        \advance\dimen7by-.5\proofrulebreadth
        \setbox0=\hbox{%
                \kern\shortenproofleft
                \ifdoubleproof
                \then   \hbox to\dimen0{%
                        $\mathsurround0pt\mathord=\mkern-6mu%
                        \cleaders\hbox{$\mkern-2mu=\mkern-2mu$}\hfill
                        \mkern-6mu\mathord=$}%
                \else   \vrule height\dimen6 depth-\dimen7 width\dimen0
                \fi
                \unhbox\proofrulename}%
        \ht0=\dimen6 \dp0=-\dimen7
\fi
\let\doll\relax
\ifwasinsideprooftree
\then   \let\VBOX\vbox
\else   \ifmmode\else$\let\doll=$\fi
        \let\VBOX\vcenter
\fi
\VBOX   {\baselineskip\proofrulebaseline \lineskip.2ex
        \expandafter\lineskiplimit\ifproofdots0ex\else-0.6ex\fi
        \hbox   spread\dimen5   {\hfi\unhbox\proofabove\hfi}%
        \hbox{\box0}%
        \hbox   {\kern\dimen2 \box\proofbelow}}\doll%
\global\dimen2=\dimen2
\global\dimen3=\dimen3
\egroup 
\ifonleftofproofrule
\then   \shortenproofleft=\dimen2
\fi
\shortenproofright=\dimen3
\onleftofproofrulefalse
\ifinsideprooftree
\then   \hskip.5em plus 1fil \penalty2
\fi
}

\newtheorem{theorem}{Theorem}[section]

\newenvironment{proofof}[1]{\begin{trivlist}\item[\hskip\labelsep{\sc
         Proof~of~{#1}.\ }]}{\hspace*{\fill} {\sc qed}\end{trivlist}}

\newcommand{\ilmo}{\sf{ILM_0}}
\newcommand{\ilm}{\sf{ILM}}
\newcommand{\ilvan}[1]{\mathfrak{Il}(#1)}
\newcommand{\ilg}{\ilvan{{\sf all}}}
\newcommand{\il}{\sf{IL}}
\newcommand{\lol}{\sf{L}}
\newcommand{\ilw}{\sf{ILW}}
\newcommand{\ilpo}{{\sf ILP}_0}
\newcommand{\ilpow}{{\sf ILP_0W^\ast}}
\newcommand{\ilp}{{\sf ILP}}

\newcommand{\tupel}[1]{\langle #1 \rangle}

\newcommand{\et}{\,\mbox{\scriptsize $\wedge$}\,}
\newcommand{\vel}{\,\mbox{\scriptsize $\vee$}\,}

\newcommand{\lmorph}[1]{\stackrel{#1}{\longrightarrow}}
\newcommand{\lver}{\sf {IL_{set}}}
\begin{document}
\title{The Interpretability Logic\\
 Of {\em All}\/ Reasonable Arithmetical Theories\\
The New Conjecture}
\author{Joost J. Joosten and Albert Visser}
\date{2000}

\maketitle

\begin{abstract}
\noindent This paper from 2000 is a presentation of a status qu{\ae}stionis at that tiime, to wit
of the problem
of the interpretability logic of {\em all}\/ reasonable arithmetical theories.
We present both the arithmetical side and the modal side of the question.
\end{abstract}

\section{Introduction}

What challenges does the future have in store for us? When talking
provability and interpretability
logic, we are in the happy position of being able to give a pretty definite
answer.
Three great problems is what we are facing. The first ---studied by R.\
Verbrugge and
A.\ Berarducci, see \cite{verb:effi93}, \cite{bera:prov93}--- is the
problem of the
provability logics of Buss' ${\sf S}^1_2$\/ and Wilkie \& Paris'
$I\Delta_0+\Omega_1$.
 The second
 is the problem of the provabilitity logic of Heyting's Arithmetic
---studied by A.\ Visser and R.\ Iemhoff, see \cite{viss:eval85},
\cite{viss:prop94}, \cite{viss:nnil95}, \cite{viss:rule99},
\cite{iemh:moda01}).
The third problem, the problem explained in this paper, is the problem of
the interpretability
logic of all reasonable arithmetical theories.

In this article, the current status of the problem will be presented. The
paper provides the necessary
definitions and a detailed explanation of the latest conjecture. It will be
made evident
that the problem is a good problem in that it intertwines modal and
arithmetical ideas.

We did our best to make this exposition
accessible to all readers with a modicum of mathematical sophistication.
The next subsection is a brief introduction to interpretations.

\subsection{What is an Interpretation?}

The interpretations we are interested in are {\em relative
interpretations}\/ in the sense
of Tarski, Mostwoski and Robinson (see \cite{tars:unde53}). Consider theories
$U$\/ with language ${\cal L}_U$\/ and $T$\/ with language ${\cal L}_T$.
For the moment we
assume that ${\cal L}_U$\/ is a relational language. An interpretation
$\cal K$\/ of
$U$\/ in $T$\/ is given by a pair $\tupel{\delta(x),F}$. Here $\delta(x)$\/
is an ${\cal L}_T$-formula
representing the {\em domain}\/ of the interpretation.\footnote{More
generally, we can
use $\delta(\vec x)$, using several variables to represent one object. We are
mainly interested in theories with sequence coding in which we can restrict
ourselves to $\delta$\/
with just one free variable.}
 $F$\/ is  a mapping that associates to each relation symbol $R$\/ of
${\cal L}_U$\/ with arity $n$\/ an ${\cal L}_T$-formula
$F(R)(x_1,\cdots,x_n)$. Here $x_1,\ldots,x_n$\/
are suitably chosen free variables. We translate the formulas of ${\cal
L}_U$\/ to the
formulas of ${\cal L}_T$\/ as follows:
\begin{itemize}
\item ${\cal K}(R(y_1,\cdots,y_n)):= F(R)(y_1,\cdots,y_n)$,\\
(We do not demand that identity is translated as identity.)
\item ${\cal K}$\/ commutes with the propositional connectives,
\item ${\cal K}(\forall y\,A) := \forall y\,(\delta(y)\to {\cal K}(A))$,
\item ${\cal K}(\exists y\,A) := \exists y\,(\delta(y)\et {\cal K}(A))$,
\end{itemize}
\noindent There are some trifling details ---e.g.\ about avoiding variable
clashes---
 that we ignore here. In case ${\cal L}_U$\/ contains functionsymbols, we first
apply the usual algorithm to eliminate functionsymbols to translate ${\cal
L}_U$\/
to a corresponding relational language and {\em then}\/ we apply the
translation sketched above.
(For an attempt to get all the details right, see \cite{viss:over98}.)
Finally, we
demand of interpretations  that for all sentences $A$\/ which are
universal closures of axioms of $U$, we have $T\vdash{\cal K}(A)$.

We will write ${\cal K}:T\rhd U$\/ for {\em $\cal K$\/ is an interpretation
of $U$\/
in $T$}. An alternative notation, which is more suitable if
we want to study the category of interpretations, is $U\lmorph{\cal K}T$.
We write $T\rhd U$, for {\em ${\cal K}:T\rhd U$, for some ${\cal K}$}.

Interpretations are used for various purposes: to prove relative
consistency, conservation results
and undecidability results. The syntactical character of interpretations
has the obvious advantage
that it allows us to convert proofs of the interpreted theory in an
efficient way
into proofs of the interpreting theory.
Examples of relative interpretations are e.g.\ the interpretation
of arithmetic in set-theory, the interpretation of elementary syntax in
arithmetic,
the interpretation of ${\sf PA}+{\sf incon}({\sf PA})$ in {\sf PA}.

Let's forget, for a brief moment, about interpretations. Let's think about
e.g.\ the construction of a model of two dimensional elliptic space in
Euclidean three
dimensional space. This is a construction inside the standard model of
Euclidean geometry,
which is modulo isomorphism the unique model of the second order version of
three dimensional
Euclidean geometry, of (modulo isomorphism) the standard model of two
dimensional
elliptic geometry. We construct this model by stipulating that `point' in
the new sense
will be {\em line through a given point}\footnote{To make this work in our
set-up,
we have to assume a version of Euclidean geometry with a constant for a
point and
an axiom stating that the point is indeed a point. To avoid the necessity
of such
inelegant stipulations we have to improve a bit on
our present definition of interpretation.}, `line' in the new sense is {\em
plane through the
given point}, `incidence of point and line' is {\em the line representing
the point
is in the plane representing the line}, etc. If you inspect the
construction, you will
see that it just uses the resources of the first order theory of three
dimensional
Euclidean geometry. Thus it provides a uniform way of transforming models of
(first order) three dimensional Euclidean geometry into models of (first order)
two dimensional elliptic geometry. A still closer inspection shows that our
construction
can be viewed as a purely syntactical transformation. It provides a
relative interpretation of
(first order) two dimensional elliptic geometry in (first order) three
dimensional Euclidean
geometry.

We can capture the relation of interpretations and model
constructions as follows. Let $\mathfrak{Mod}(V)$\/ be the class of models
of a theory $V$.
An interpretation $\cal K$\/ of $U$\/ in $T$\/ provides a uniform way to
build internal models of
$U$\/ inside models of $T$.
Thus $\cal K$\/ provides us with a function, say $\mathfrak{Mod}({\cal
K})$, from
$\mathfrak{Mod}(T)$ to $\mathfrak{Mod}(U)$. Thus defined $\mathfrak{Mod}$
is a contravariant
functor from $\mathfrak{Theory}$, the category of theories and
interpretations, to
$\mathfrak{Class}$, the category of definable classes and definable
functions between classes.
The idea of an interpretation as an `internal model given in a uniform way'
is an important
heuristic in thinking about interpretations: the mind craves reality and
visualisation
 rather than syntax. We will exploit this heuristic in what follows. E.g.\
we will speak about
one interpretation $\cal K$\/ being an end-extension of another one $\cal
M$, meaning that
in every model the internal model associated to $\cal K$\/ is an
end-extension of the
internal model of $\cal M$\/ in a uniform way.

\subsection{So What's Reasonable?}

We will be interested in interpretability of reasonable arithmetical
theories. More specifically
we will be interested in what such theories have to say about
interpretability in these theories
themselves.
So what are reasonable arithmetical theories?

A theory, for our present purposes,
 is a predicate logical theory axiomatized by axioms in an axiom set that
is given by
arithmetical formula $\alpha$. Unless stated otherwise we assume that
$\alpha$\/ is simple,
say a predicate corresponding with a class that is decidable by a p-time
algorithm. Note that
our specification makes
 {\em theory}\/ an intensional notion, since we consider theories also from
the point of view
of theories: $\alpha$\/ and $\beta$\/ may specify the same axiom set, but a
theory $U$\/
thinking about $T_\alpha$\/ and $T_\beta$\/ need not be aware of that.

`Arithmetical' has a primary and a secondary meaning. In the primary
meaning an arithmetical theory
is an extension of Robinson's Arithmetic {\sf Q} in the usual language of
arithmetic with
$\underline 0$, $S$, $+$ and $\times$. (We will often use $\cdot$ instead
of $\times$.) In the
secondary meaning, an arithmetical theory is a pair $\tupel{T,{\cal N}}$, where
${\cal N}:T\rhd {\sf Q}$. In other words an arithmetical theory is a theory
with designated predicates
representing the natural numbers, representing zero, etc. e.g. $\tupel{{\sf
ZF},{\cal N}}$,
where $\cal N$\/ is the usual set-theoretical representation of the natural
numbers, is an
arithmetical theory. The intended meaning of arithmetical theory
in this paper is the secondary one.

The caution concerning the explicit designation of the natural numbers
is necessary, since not all interpretations of number theory are provably
isomorphic in a given theory.
Thus the following three statements are equally true.
\begin{enumerate}
\item ${\sf con}({\sf ZF})$ is independent of {\sf GB} (= G\"odel-Bernays
Set Theory).
\item {\sf GB} proves ${\sf con}({\sf ZF})$.
\item {\sf GB} proves $\neg\,{\sf con}({\sf ZF})$.
\end{enumerate}

\noindent Here ${\sf con}({\sf ZF})$ abbreviates a fixed arithmetical
sentence, but we vary, in the examples, the designated
set of natural numbers. The correct formulation of our statements is:
\begin{enumerate}
\item ${\sf con}({\sf ZF})$ is independent of $\tupel{{\sf GB},{\cal N}}$.
\item $\tupel{{\sf GB},{\cal I}}$ proves ${\sf con}({\sf ZF})$.
\item $\tupel{{\sf GB},{\cal K}}$ proves $\neg\,{\sf con}({\sf ZF})$.
\end{enumerate}

\noindent Here $\cal N$\/ is the usual interpretation of the natural
numbers in {\sf ZF} lifted to {\sf GB}.
$\cal I$\/ is a definable cut of the $\cal N$-numbers and $\cal K$\/ is a
suitable interpretation
built using a syntactic variant of the Henkin construction. Both $\cal N$\/
and $\cal I$\/ are standard
in that they represent the ordinary natural numbers (modulo isomorphism) in
the standard model.
$\cal K$, of course, cannot represent the standard numbers inside any model.

Rather than viewing the
possibility of having different sets of numbers as a nuisance, we will make
grateful use of it
by switching between choices of what `the numbers' are.

There is some arbitrariness in our singling out arithmetic
as the thing we are interested in, especially since representations of
syntax play
such a large role in G\"odelean metamathematics. We could as have well
decided to speak about
{\em syntactical}\/ theories, counting e.g.\ $\tupel{{\sf PA},{\cal S}}$
as a syntactical theory, where ${\cal S}$\/ is a
designated interpretation of some reasonable theory of elementary syntax.
There are two
reasons we will make the traditional choice to speak about arithmetic:
first simply because
it is the traditional choice ---changing it will cause confusion---,
secondly the methodology
of definable cuts is easier to understand in the context of arithmetic.

What is reasonable? It means {\em at least}: strong enough to verify the
minimal principles we are interested in. Take, e.g.\ the principle
that tells us
that if something is provable, then it's provable that it's provable. To
verify the principle
in the obvious way, we need $\Delta_0$-induction, plus the totality of the
function
$\omega_1$, where $\omega_1(x)=2^{(\log_2x)^2}$. This principle is called
$\Omega_1$.
Another principle is the one stating that interpretations can be composed,
i.e.\
if ${\cal K}:U\rhd V$\/ and ${\cal M}:V\rhd W$, then ${\cal M}\circ{\cal
K}:U\rhd W$.
To verify it we need
$\Sigma$-collection, also known as $B\Sigma_1$,
the principle $\forall x\leq a\,\exists y\; A \to \exists b\,\forall x\leq a\,
\exists y\leq b\;A$, where $A$\/ is $\Sigma_1$.

Thus we demand that reasonable arithmetical theories contain a minimal
arithmetical theory
{\sf Basic}. Formally: a reasonable arithmetic is of the form
$\tupel{T,{\cal N}}$, where
${\cal N}:T\rhd{\sf Basic}$. The most plausible choice for {\sf Basic} at
the time of writing is
$I\Delta_0+\Omega_1+B\Sigma_1$. (See \cite{viss:form91}, for some shameless
trickery to get rid of
the assumption $B\Sigma_1$.)

The second demand that we pose, has to do with the coherence of the
theories. Given two theories
$T,U$, we could take a disjoint union $T\oplus U$\/ in such a way that the
$T$- and the $U$-objects
have no recognizable interaction at all. So if our numbers are `confined
in' $T$, they will not be able to
`interact' in any way with the $U$-part. What we demand is that our theory
is sequential: it
should contain (in the sense of interpretability) a theory of sequences of
{\em all}\/ objects the
theory can talk about.
Here the lengths ${\sf length}(\sigma)$ and projections
$(\sigma)_i$\/ are taken from an initial segment of the designated numbers.
Sequentiality is important to make e.g.\ the construction of partial
truthpredicates possible.
For more on sequentiality, see e.g.\ \cite{haje:meta91}.

The third demand is not really a demand but a programmatic point. We should
keep the
answer to the question what a reasonable theory is, somewhat indefinite.
The class of reasonable
theories is that class of theories that allows a beautiful answer to the
question what the interpretability logic of all reasonable theories is.
E.g.\ it could happen that only the theories that contain the axiom {\em
that exponentiation is total}\/
have a nice logic. Well, in that case we say that those theories are the
reasonable ones.

\subsection{Approaches to Interpretability}

What could the metamathematical study of interpretability and
interpretations look like?
One idea is to study {\em degrees of interpretability}. Interpretability
yields a partial
preorder on theories. Dividing the associated equivalence relation out we
get a degree-theory.
Degree-theory has been studied by P.\ Lindstr\"om and C.\ Bennet (see e.g.\
\cite{lind:aspe97}
and \cite{benn:orde86}) and by  V.\ {\v{S}}vejdar (see \cite{svej:degr78}).
The work on degrees was very fruitful as a generator of methods and techniques.
Some of these techniques have been adapted for application in
interpretability logic.

We feel that it could be very fruitful to extend the degree-theoretic
approach to the
study of suitable categories of interpretations. The more expressive
category-theoretical language
might be better suited to express certain basic insights concerning
interpretability.
There were some attempts to initiate such a study, but these attempts did
yield less
than satisfactory results. Some further experimentation is needed to
isolate the right
categories.

The approach to interpretability that is the focus of this paper is the
modal study of
interpretability.\footnote{Lev Beklemishev places this kind
of study between structural prooftheory which studies specific proof
systems and proofs,
and recursion theoretic prooftheory where theories are considered as RE
sets of theorems.
Here we abstract away from many details of the proof system and from
detailed proofs, however
e.g.\ the defining formula of the set of axioms is a feature that can make
a difference. Perhaps one could say that
the modal study of provability and interpretability is part of {\em
intensional prooftheory}.}
The modal language has the advantage of expressiveness, but there are
costs. First modal logic is about `propositions' not about theories. This
means that
we cannot directly study the relations ${\cal K}:T\rhd U$\/ or even $T\rhd
U$. What we study is
the relation $A\rhd_T B$, which is defined as follows:
 \begin{itemize}\item $A\rhd_TB\;\; :\Leftrightarrow
(T+A)\rhd(T+B)$.\end{itemize}
Here $T$\/ is the base theory. We speak of (sentential) interpretability
{\em over}\/ $T$.
Secondly, we are interested in iterating the modal connectives. We want to
allow things like $(A\rhd B)\rhd C$. This means that our research is
restricted to base theories
$T$\/ that have sufficient coding ability to formalize a decent amount of
reasoning concerning
interpretations. This restriction is substantial since lots of important
interpretations
fall outside  the scope of our investigation. If we pay the costs, there
are some gains.
\begin{enumerate}
\item We have a modal language that is rich enough to articulate both
       the incompleteness theorems and the model existence lemma, which is the
       heart of the completeness theorem.
\item Some substantial reasoning concerning interpretability can be
executed in this
      modal logic.
\item The Kripke model theory of the logic is highly interesting qua modal
logic.
\item The arithmetical side of the study involves substantial arithmetical
insights.
        As we will see, in an indirect way, our logic can talk about large
and small numbers.
\end{enumerate}
Before we introduce the modal logics, we interpolate a brief introduction
to some salient arithmetical
facts.

\section{Parvulae Arithmeticae}

\subsection{Coding}

Since the function $\omega_1=\lambda x.2^{(\log_2x)^2}$
is present in our basic system of arithmetic {\sf Basic},
we have p-time computable functions available. Having these,
arithmetization of syntax becomes
a piece of cake. The most obvious g\"odelnumbering of strings in a given
alphabet is
also the best one. We enumerate first the strings of length 0, then the
strings of length 1,
and so on. The strings of the same length we order alphabetically. We
assign to each string
as g\"odelnumber simply its ordernumber in the sequence so obtained.
It turns out, using a trick due to Smullyan, that operations on strings
like concatenation
can be easily arithmetically represented under this coding. An important
insight is the
elementary fact that the g\"odelnumber of a string $\sigma$\/ is of order
$A^{{\sf length}(\sigma)}$,
where $A$\/ is the cardinality of the alphabet. We will code syntactical
elements, like formulas
and proofs, by writing them out and then taking the code of the resulting
string.

We will write $\Box_TA$\/ for the arithmetization of {\em $T$\/ proves
$A$}\/ and $A\rhd_TB$\/
for the arithmetization of {\em $T+A$\/ interprets $T+B$}. If $A$\/
contains a free variable $x$,
$\Box_TA$\/ is the arithmetization of {\em the result of substituting the
numeral of $x$\/ in $A$
for ``$x$'' is provable in $T$'}. Further conventions are similar.

\subsection{Efficient Numerals}

It is definitely
not a good idea to represent the number $n$\/ by the numeral
$\overbrace{S\cdots S}^n\underline 0$. The g\"odelnumber of this numeral
will be
of order $2^{cn}$, for a fixed constant $c$. So the function sending a number
to the code of its numeral will be exponential. Exponentiation is not
generally available in
{\sf Basic}. Hence we will use binary numerals instead. These are defined by
${\sf num}(0):= \underline 0$,
${\sf num}(2n+1):= S(SS\underline 0\cdot{\sf num}(n))$,
${\sf num}(2n+2):= SS\underline 0\cdot{\sf num}(n+1)$. This representation
has the happy consequence that the g\"odelnumber of the numeral of $n$\/ is of
order $2^{c\log_2n}$, i.e.\ $n^k$, for some fixed standard $k$.

\subsection{Numbers Large and Small}

We have to face the basic fact that we are going to use theories that do not
have full induction. Note that also quite strong theories may lack full
induction, e.g.\
$\tupel{{\sf GB},{\cal N}}$.
Thus, in our theories, it may happen that we have definable sets of numbers
containing 0
and closed under successor such that the theory doesn't think this set
contains all (designated) numbers.
In some cases the theory will even positively know this set does not
contain all numbers.
Such  definable sets of numbers, closed under successor but not necessarily
containing all designated numbers, play an important role in the
metamathematical study of arithmetics.
For many purposes it is convenient to put stronger demands on these sets:
we ask that they
are {\em definable cuts}. Let $\tupel{T,{\cal N}}$\/ be an arithmetical theory.
Here ${\cal N}=\tupel{\delta,F}$.

\bigskip\noindent
An ${\cal L}_T$-formula $I$\/ is/presents a $\tupel{T,\cal N}$-cut iff
$T$\/ proves that:
\begin{enumerate}
\item $Ix \to \delta x$,
\item $(Ix \et {\cal N}(x=y)) \to Iy$,
\item $I$\/ is downwards closed under $<$, i.e.\\
    $(I(x) \et {\cal N}(y<x)) \to I(y)$,
\item $I$\/ is closed under $0$, $S$, $+$, $\times$ and $\omega_1$, i.e.
\begin{enumerate}
\item ${\cal N}(x=0)\to Ix$,
\item $(Ix\et{\cal N}(Sx=y))\to Iy$,
\item $(Ix\et Iy \et {\cal N}(x+y=z)) \to Iz$,
\item  $(Ix\et Iy \et {\cal N}(x\cdot y=z)) \to Iz$,
\item $(Ix\et{\cal N}(\omega_1x=y))\to Iy$.
\end{enumerate}
      Note that `$\omega_1x=y$' is, in the usual set-up, an abbreviation of
a complex formula.
\end{enumerate}
We will sometimes write `$x\in I$' for `$Ix$'.

Using a wonderful trick invented by Solovay \cite{solo:inte76}, we can
always `shorten' a definable
set of numbers, $T$-provably closed under successor to a $T$-cut.
Cuts can be considered as `notions of smallness': the numbers inside the
cut are `small', the ones above it `big'.

We will consider cuts themselves as interpretations of arithmetic,
confusing the
cut $I$\/ with the interpretation $\tupel{I,F}$, where $F$\/ is the
interpretation function associated
with $\cal N$. It is easy to see that $\tupel{I,F}$ is indeed an
interpretation.

\bigskip\noindent
{\em From this point on, we will often suppress the designated cut $\cal
N$, writing as if
$\cal N$\/ were the identity interpretation.}

\bigskip\noindent
A startling fact about cuts is the {\em outside big, inside small}\/ principle.
Even if $T$\/ may fail to believe that every number is in the $T$-cut $I$, we
do have:

\begin{theorem}\label{oubiinsm} $T\vdash\forall x\,\Box_T\;x\in I$.
\end{theorem}

\noindent The point is that we can have big proofs showing that big numbers
are small.
Here is a somewhat more elaborate proofsketch.

\begin{proofof}{\ref{oubiinsm}}
We reason informally in $T$.
Let $\sigma$\/ be a (standard) proof of $\forall x\;(x\in I\to
S(SS\underline 0\cdot x)\in I)$.
We convert a proof $\pi$\/ of $\underline n\in I$\/ into a proof of
$S(SS\underline 0\cdot\underline n) \in I$\/ as follows.

\[ \prooftree         \pi \;\;\;\;\;
         \[ \sigma
     \justifies \underline n\in I \to S(SS\underline 0\cdot\underline n)\in I
      \thickness=0.08em
  \using \forall E \]
  \justifies S(SS\underline 0\cdot\underline n) \in I
  \thickness=0.08em
  \using \to E
\endprooftree \]

\noindent Similarly we convert a proof of $\underline n\in I$
into a proof of $(SS\underline 0\cdot\underline n) \in I$.
 Clearly a proof of $\underline n\in I$\/ will use in the order of
$\log_2 n$\/ steps. The number of symbols in a step of the proof
can be estimated by $a\log_2 n +b$\/ for fixed standard numbers $a$\/ and $b$.
The number of symbols in the proof will be estimated by: $\log_2 n\cdot (a
\log_2 n +b)$.
Let $c:=a+b$. We can replace our estimate by: $c\cdot(\log_2(n+2))^2$.
So the size of the G\"odelnumber of the proof will be estimated by
$2^{c\cdot(\log_2(n+2))^2}=(\omega_1(n+2))^c$. The function
$\lambda x.(\omega_1(x+2))^c$ is present in {\sf Basic}.\end{proofof}

\subsection{Cuts and Interpretations}\label{cutint}

If we think of an interpretation $\cal K$\/ of $U$\/ in $T$\/ as an inner
model of $U$\/
inside a model of $T$, we can ask how the $T$-numbers do compare to the
$U$-numbers as
seen via $\cal K$.
To be pedantically precise, if the $T$-numbers are given by $\cal N$\/ and
if the $U$-numbers are
given by $\cal M$, how does the internal model of {\sf Basic} given by
$\cal N$ compare to the
internal model of {\sf Basic} given by ${\cal K}\circ{\cal M}$?

In case $T$\/ has full induction, the answer is simple. Let's for the moment
step into the outside world and remind ourselves of a basic fact about
non-standard
models of arithmetic. The natural numbers form (modulo embedding) an
initial fragment
of every non-standard model. I.o.w.\
every non-standard model is an end-extension of the standard model.
We can prove that fact by defining the embedding of the natural numbers
into the
non-standard model by external recursion
  and by subsequently proving the desired properties of the embedding by
external induction.
Essentially the same argument can be repeated in $T$. We define the embedding,
using the fact that we are supposed to have sequences of objects, by an
explicit $T$-formula
and verify its properties with $T$-induction. The upshot is that the
$U$-numbers as seen
via $\cal K$\/ form an end-extension of the $T$-numbers.

Now what happens if $T$\/ doesn't have full induction? The answer to this
question
has been provided by Pavel Pudl\'ak in his fundamental paper
\cite{pudl:cuts85}. Here it is.
There is a $T$-cut $I$\/ such that the $U$-numbers as viewed via $\cal K$\/
are an
end-extension of $I$.

\bigskip\noindent
The argument for Pudl\'ak's theorem is a refinement of the usual argument
sketched above: {\em where we lack induction, we compensate by switching to
a smaller cut}.

\bigskip\noindent {\bf A Closer Look}\\
To convince the reader that the statement of Pudl\'ak's theorem makes sense, we
spell out the result in the pedantic mode. Remember that we assumed that
$T$\/ was really
$\tupel{T,{\cal N}}$\/ and $U$\/ was really $\tupel{U,{\cal M}}$. Now
${\cal Q}:={\cal K}\circ{\cal M}$\/
 is an interpretation of {\sf Basic} in $T$, representing the $U$-numbers
as viewed by $T$\/ via $\cal K$.
Pudl\'ak's Theorem tells us that there is a $\tupel{T,{\cal N}}$-cut $I$\/
such that there is, verifiably in $T$, a definable embedding of $I$ into an
initial segment
of  the $U$-numbers as viewed in $T$\/ via $\cal Q$.
This means that there is a $T$-formula $E$\/ such that $T$\/ proves:
\begin{enumerate}
\item
$Exx' \to (Ix \et \delta_{\cal Q}(x'))$ \\($E$\/ is a relation between
$I$\/ and $\delta_{\cal Q}$),
\item
$(Exx' \et {\cal N}(x=y) \et {\cal Q}(x'=y')) \to Eyy'$ \\
($E$\/ is a congruence w.r.t.\ the relevant `identities'),
\item
$Ix \to \exists x'\;Exx'$\\ ($E$\/ is total on $I$),
\item
$(Exx' \et Exy') \to {\cal Q}(x'=y')$ \\($E$\/ is a function),
\item
$(Exx' \et Eyx') \to {\cal N}(x=y)$\\ ($E$\/ is injective),
\item
$(Exx'\et {\cal Q}(y'<x')) \to \exists y \;({\cal N}(y<x) \et yEy')$\\ (The
$E$-image of $I$\/ is
downwards closed in $\cal Q$),
\item
$(Exx' \et {\cal N}(x=0) \to {\cal Q}(x'=0)$\\ ($E$\/ commutes with 0),
\item
$(Exx'\et Eyy' \et {\cal N}(Sx=y)) \to {\cal Q}(Sx'=y')$\\ ($E$\/ commutes
with $S$),
\item
$(Exx'\et Eyy' \et Ezz' \et {\cal N}(x+y=z)) \to {\cal Q}(x'+y'=z')$\\
($E$\/ commutes with $+$),
\item
$(Exx'\et Eyy' \et Ezz' \et {\cal N}(x\cdot y=z)) \to {\cal Q}(x'\cdot
y'=z')$ \\
($E$\/ commutes with $\times$).
\end{enumerate}
The image, say $J$, of the $\tupel{T,{\cal N}}$-cut $I$\/ is easily seen to
be a $\tupel{T,{\cal Q}}$-cut.
$T$\/ shows that $E$\/ is an isomorphism between $I$\/ and $J$. $J$\/
is  not generally internally definable in $\cal K$,
i.o.w.\ there need not be an ${\cal L}_U$-formula $G$\/ such that
$T\vdash \exists x\,Exx' \leftrightarrow {\cal K}(Gx')$.

\section{Interpretability Logic Explained}

\subsection{Description of the System {\sf IL}}

The language of interpretability logic, ${\cal L}_{\sf int}$,
is the language of modal propositional logic extended
with a binary modal operator $\rhd$. We read
$A\rhd B$\/ as: $A$\/ interprets $B$.
We will write $\Diamond A$ as an abbreviation of $\neg\Box\neg$.

Let $U$\/ be a reasonable arithmetical theory. An interpretation $(.)^\ast$ of
${\cal L}_{\sf int}$ into $U$\/ maps the atoms on sentences of ${\cal L}_U$,
commutes with the propositional connectives and satisfies: \begin{center}
$(\Box A )^\ast :=\Box _U A^\ast $ and
$(A\rhd B )^\ast :=A^\ast \rhd_U B^\ast$.\end{center}

\noindent
We study the interpretability principles valid in theories $U$, asking
ourselves for which $C$\/
in the modal language
 we have $U\vdash C^\ast$, for all $(.)^\ast$ and asking ourselves
which principles are valid in all reasonable theories. The set of
principles valid in $U$\/ is called
$\ilvan{U}$. The set of principles valid in {\em all}\/ reasonable theories
will be called
$\ilvan{{\sf all}}$.\footnote{For the modal language restricted to the {\em
unary}\/
 connective $\top\rhd A$\/
in combination with $\Box$, the problem of the interpretability logic of
all theories has been
solved by Maarten de Rijke, see his \cite{rijk:unar92}.}

We introduce the basic modal logic
{\sf IL}. The principles of this logic are arithmetically sound for a wide
class of theories and
for various interpretations of its main connective $\rhd$.\footnote{We can
also interpret
$\rhd$ as partial conservativity w.r.t.\ a suitable class of formulas.} The
theory is
arithmetically incomplete for all known arithmetical interpretations.
The motivation for studying this specific set of axioms
comes from its {\em modal}\/ simplicity and elegance.

{\sf IL} is the smallest logic in ${\cal L}_{\sf int}$ containing the
tautologies of propositional logic,
closed under modus ponens and the following rules. (A principle is just a
rule with empty antecedent.)
\begin{itemize}
\item[{\sf L1}]	$\vdash A \;\Rightarrow\; \vdash \Box A  $
\item[{\sf L2}]	$\vdash \Box (A  \to B  ) \to (\Box A  \to \Box B  )$
\item[{\sf L3}]	$\vdash \Box A  \to \Box \Box A  $
\item[{\sf L4}]	$\vdash \Box (\Box A  \to A  )\to \Box A  $
\item[{\sf J1}]	$\vdash \Box (A  \to B  ) \to A  \rhd B  $
\item[{\sf J2}]	$\vdash (A  \rhd B  \;\et\; B  \rhd C ) \to A  \rhd C $
\item[{\sf J3}]	$\vdash (A  \rhd C  \;\et\; B  \rhd C  ) \to (A  \vel B
)\rhd C  $
\item[{\sf J4}]	$\vdash A  \rhd B  \to (\Diamond A  \to \Diamond B  )$
\item[{\sf J5}]	$\vdash \Diamond A  \rhd A $
\end{itemize}
{\sf L1-4} are the well-known principles of L\"ob's Logic.
{\sf IL} is certainly valid in all reasonable theories $U$. We will provide the
arithmetical justifications of the principles in subsection~\ref{checkil}.

D.\ de Jongh and A.\ Visser proved that {\sf IL} has {\em unique and
explicit fixed points}.
See \cite{dejo:expl91}.
No characterization of the {\em closed fragment} of {\sf IL} has been given.
{\sf IL} satisfies {\em interpolation}, see \cite{arec:inte98}.
De Jongh and Veltman prove a {\em modal completeness theorem}\/ w.r.t.\
Veltman models. See \cite{dejo:prov90}.

Here is a sample of {\sf IL}-reasoning. We prove: $\vdash
A\rhd(A\et\Box\neg \,A)$.
First, by {\sf L1-4}, we can derive, taking the contraposition of
{\sf L4}: $\vdash \Diamond A \to \Diamond (A\et\Box\neg\,A)$. So, by {\sf L1}
and {\sf J1}, we find: $\vdash \Diamond A\rhd \Diamond (A\et\Box\neg\,A)$.
Applying {\sf J5} and {\sf J2}, we get: (a) $\vdash \Diamond A\rhd(A\et
\Box\neg\,A)$.
We also have, by {\sf L1} and {\sf J1}: (b) $\vdash A\rhd
((A\et\Box\neg\,A)\vel\Diamond A)$
and (c): $\vdash(A\et\Box\neg\,A)\rhd(A\et\Box\neg\,A)$. Applying {\sf J3}
and {\sf J2} to (a), (b) and (c) we arrive at the desired result.

Putting $\top$\/ for $A$\/ in the principle we just derived, we see that it
follows that
one can construct, in a uniform way,
 inside every model of a given arithmetical theory $T$\/ an internal
model of $T+{\sf incon}(T)$.

\subsection{The Arithmetical Validity of {\sf IL}}\label{checkil}

\subsubsection*{Verification of the {\sf L}-principles}

It is well known that the principles of L\"ob's Logic can be derived in Buss'
${\sf S}^1_2$ (see \cite{buss:boun86}) or in Wilkie \& Paris'
$I\Delta_0+\Omega_1$
(see \cite{wilk:sche87}). Since {\sf Basic} is
supposed to extend ${\sf S}^1_2$, we are done. The proof of
{\sf L3} is by induction on the subformulas of $\Box_TA$, using the fact that
$\Box_TA$\/ is a $\exists\Delta^b_0$-predicate. A $\Delta^b_0$-formula only
has logarithmically
bounded quantifiers.

It is a remarkable fact that {\sf L3} is doubly redundant in {\sf IL}. By a
clever argument, due to Dick de Jongh, we can derive {\sf L3} from
{\sf L1,2,4}. However this redundancy is not arithmetically helpful,
since the usual arithmetical verification of L\"ob's axiom {\sf L4} uses the
validity of {\sf L3}.\footnote{We can derive {\sf L4} without using
{\sf L3} by employing a surprising argument of Kreisel (presented in
\cite{smor:inco77}).
However, this argument includes the verification of {\sf J5}.} The second way
is to derive {\sf L3} from {\sf J5} and {\sf J4}. The striking thing about
this alternative proof
is that it provides a really different way to obtain a $T$-proof of
$\Box_TA$\/ from a $T$-proof
of $A$.\footnote{Andr\'e van Kooy showed in his masters thesis (Department of
Philosophy, Utr\-echt University) that for finitely axiomatized
theories in a relational language one can make the
transformation of a proof of $A$\/ into a proof of the provability of $A$
{\em linear time}.
This seems to be only possible via the {\sf J4,5}-route.}
{\sf L3} is one of those cases where we have one fact but two insights.

Both styles of proofs of {\sf L3}, yield on inspection sharper results, like:
\begin{itemize}
\item $T\vdash \Box_TA\to\Box_T\Box_T^IA$\\
Here $I$\/ is any $T$-cut. We write $\Box_T^IA$ for $I(\Box_TA)$ ---note
that we need only to relativize the unbounded existential quantifier of
$\Box_TA$\/ to $I$.)
\item $T\vdash \Box_TA\to\Box_T\triangle_TA$\\
Here $\triangle$ stands for either cutfree, Herband or tableaux
provability.\end{itemize}

\noindent The derivation of the first strengthening in the {\em induction
on subformulas}\/
style, runs as follows.

We obtain at a certain point
$\Box_T{\sf proof}_T(x,\underline{{\sf gn}(A)})$. (Here {\sf gn} is the
g\"od\-el\-numb\-er\-ing function.)
The {\em outside big, inside small}\/ principle
tells us that $\Box_Tx\in I$. Ergo, $\Box_T\Box_T^IA$.

\subsubsection*{Verification of {\sf J1}}

The validity of {\sf J1} is witnessed by the identity interpretation {\sf ID}.

\subsubsection*{Verification of {\sf J2}}

If ${\cal K}:A\rhd_TB$\/ and ${\cal M}:B\rhd_TC$, then $({\cal M}\circ{\cal
K}):A\rhd_TC$.

\subsubsection*{Verification of {\sf J3}}

{\sf J3} is valid, since,
given any two interpretations $\cal K$\/ and $\cal M$\/ and any sentence $A$,
we can construct an interpretation ${\cal K}[A]{\cal M}$, the disjoint
 $A$-sum of $\cal K$\/ and $\cal M$,
that behaves like $\cal K$\/ if $A$\/ and like $\cal M$\/ if $\neg A$.
We take:
\begin{itemize}
\item
$\delta_{{\cal K}[A]{\cal M}}(x) :=
((\delta_{\cal K}(x)\et A)\vel (\delta_{\cal M}(x)\et \neg\,A))$,
\item
$({\cal K}[A]{\cal M})(P)(\vec x) := (({\cal K}(P)(\vec x)\et A) \vel
({\cal M}(P)(\vec x)\et \neg\,A))$
\end{itemize}
We find that if ${\cal K}:A\rhd_TC$\/ and ${\cal M}:B\rhd_TC$, then
$({\cal K}[A]{\cal M}):(A\vel B)\rhd_TC$.

\subsubsection*{Verification of {\sf J4}}

{\sf J4} tells us that relative interpretability
implies relative consistency. If we would have $\Box_T\neg\,B$\/ and ${\cal
K}:A\rhd_TB$,
then it would follow that $\Box_T(A\to \bot^{\cal K})$, and hence
$\Box_T\neg\,A$.

\subsubsection*{Verification of {\sf J5}}

{\sf J5} is the {\em interpretation existence lemma}.
It is the syntactical realization of the Henkin model existence lemma.
Inspecting the usual proof of the model existence lemma, one sees that it
involves the
construction of a set of sentences describing a model. This set can as well
be viewed as describing an
interpretation. The set is constructed as a path in a binary tree. This
path is described by a
 $\Delta^0_2$-predicate. The desired properties of the set of sentences
are verified using induction. Thus the whole argument can easily be
verified in {\sf PA}.
The construction can be executed in almost any arithmetical theory by
compensating for the lack of induction by switching to definable
cuts.\footnote{We are
not quite sure who found this fact first. It might be well have been
discovered independently
by Friedman, Pudl\'ak and Solovay.}

A moment's reflection shows that the choice of the numbers in which we
execute the Henkin construction
is irrelevant. So, in particular, this set of numbers, might very well be
some $T$-cut $I$. Thus
we arrive at the following sharpening of {\sf J5}.
\begin{itemize}
\item $T\vdash \forall I\;(\Diamond^I_TA\rhd_TA)$.
\end{itemize}
Here $\Diamond^I_TA$\/ stands for $I(\Diamond_TA)$, which is $T$-equivalent to
\[\forall x\in I\;\neg\,{\sf prov}_T(x,\underline{{\sf gn}(\neg\,A)}),\]
 where {\sf gn} is the
g\"odelnumbering function.

Note that the first strengthened version of {\sf L3}, follows easily from
{\sf J4} and the
sharpened version of {\sf J5}. We have $T\vdash\forall
I\;(\Diamond_T^IB\rhd_TB)$, and, hence,
$T\vdash\forall I\;(\Diamond_T\Diamond_T^IB\to\Diamond_TB)$. Replacing
$B$\/ by $\neg\,A$,
contraposing and cleaning up spurious negations (using {\sf L1,2}), we get:
$T\vdash\forall I\;(\Box_TA\to\Box_T\Box_T^IA)$.

\subsection{Beyond {\sf IL}}

{\sf IL} is certainly arithmetically sound. However, it is not
arithmetically complete for any
reasonable
arithmetical theory $T$ and for any known interpretation of $\rhd$.

\subsubsection*{Montagna's Principle {\sf M}}

Let us first consider Peano Arithmetic, {\sf PA}.
The theory satisfies a further principle: Montagna's Principle.

\begin{itemize}
\item[{\sf M}] $\vdash  A\rhd B \to (A\et \Box C)\rhd(B\et \Box C)$
\end{itemize}

\noindent The {\sf PA}-validity of {\sf M} was known independently to
\v{S}vejdar and Lindstr\"om.
Arithmetical completeness for the system {\sf ILM}:={\sf IL}+{\sf M} was
conjectured by A.\ Visser.
It was proved independently by V.\ Shavrukov (see \cite{shav:rela88}) and
A.\ Berarducci
(see \cite{bera:inte90}).
For nice presentations of the proof see also \cite{zamb:proo92} or
\cite{japa:logi98}.
It turns out that {\sf ILM} is sound and complete for all reasonable
arithmetical theories satisfying
full induction.\footnote{In fact the class is somewhat bigger. The reader is
referred to
\cite{viss:over98} for further elaboration.} Here we just verify the
arithmetical
validity of {\sf M}.

\bigskip\noindent Let $T$\/ have full induction. We prove the stronger
principle:
\begin{itemize}
\item $\vdash  A\rhd_T B \to (A\et S)\rhd_T(B\et S)$, for $S\in\Sigma^0_1$
\end{itemize}
Reason informally in $T$. Suppose ${\cal K}:(T+A)\rhd(T+B)$. Reason in $T$.
(So we
are in $T$\/ two deep.) Suppose $A$\/ and $S$. Consider $\cal K$. We will
certainly have ${\cal K}(C)$
for each axiom $C$\/ of $T$\/ and ${\cal K}(B)$, since we have $A$. Now the
$T+B$-numbers seen
via $\cal K$\/ are an end-extension of the $T+A$-numbers, as explained in
subsection~\ref{cutint}.
Moreover, $\Sigma^0_1$-sentences are preserved by end-extensions. Hence, we
have ${\cal K}(S)$.

\subsubsection*{The Persistence Principle {\sf P}}

The persistence principle {\sf P} is the following principle:
\begin{itemize}
\item[{\sf P}] $\vdash A\rhd B \to \Box(A\rhd B)$.
\end{itemize}
The persistence principle is valid for interpretations in finitely
axiomatized reasonable
arithmetical theories $T$. The reason is simple. Let $C$\/ be the
conjunction of the axioms of
$T$. Then, to obtain our principle, is is sufficient to verify:
\[ T\vdash \exists {\cal K}\,\Box_T(A\to {\cal K}(C\et B)) \to
      \Box_T\exists {\cal K}\,\Box_T(A\to {\cal K}(C\et B)) \]
which is obviously valid by verifiable $\exists\Delta_0^b$-completeness.

Albert Visser has shown that {\sf ILP} is arithmetically complete for each
finitely
axiomatized reasonable theory that proves {\sf Superexp}, the axiom stating
that
superexponentiation is total. See \cite{viss:inte90}.
It is definitely known that {\sf ILP} is {\em not}\/ complete
for $I\Delta_0+B\Sigma_1+{\sf Exp}$, where {\sf Exp} is the axiom stating that
exponentiation is total.

\subsubsection*{The Principle {\sf W}}

The first principle that was found to be valid in all theories that
strictly extends
{\sf IL} is the principle {\sf W}. `{\sf W}' for ``weak''.
\begin{itemize}
\item[{\sf W}]
$\vdash A\rhd B\to A\rhd (B\et \Box\neg\,A)$
\end{itemize}
For some time it was conjectured
that {\sf ILW} was $\ilvan{{\sf all}}$. This conjecture was eventually refuted.
Before we turn to the next principle, we verify {\sf W}.

\bigskip\noindent
Remember our verification of the principle $\vdash A\rhd(A\et\Box\neg\,A)$.
Now all the principles used remain valid if we relativize all the
modal operators at nesting depth 1 to a $T$-cut $I$. Thus we obtain:
\begin{itemize}
\item $T\vdash \forall I\;(A\rhd_T(A\et\Box_T^I\neg\,A))$
\end{itemize}
Now reason in $T$. Suppose, for some $\cal K$, we have
${\cal K}:A\rhd_TB$. Let $J$\/ be the
$T+A$-cut such that, in $T+A$,  $J$\/ is (isomorphic to) an initial segment
of the $T+B$-numbers
viewed via $\cal K$.
If $J$\/ is not a $T$-cut, we may replace it by $J[A]{\sf ID}$, which
certainly is a
$T$-cut. (The notation $(\cdot)[\cdot](\cdot)$ was introduced when
explaining {\sf J3}.)
So we can assume that $J$\/ is a $T$-cut.
Now we have (a) ${\cal K}:(A\et\Box_T^J\neg\,A)\rhd_T(B\et\Box_T\,\neg A)$,
since $J$\/ is initial in $\cal K$\/ and since $\Sigma^0_1$-sentences are
upwards persistent.
By our previous consideration we have (b) $A\rhd_T(A\et\Box_T^J\neg\,A)$.
Composing (a) and (b) we
arrive at the desired result.

\subsubsection*{The Principle ${\sf M}_0$}

The next principle that was discovered is the principle ${\sf M_0}$.
\begin{itemize}
\item[${\sf M_0}$] $\vdash A\rhd B \to (\Diamond A\et \Box C)\rhd (B\et\Box C)$
\end{itemize}

\noindent Here is the arithmetical verification. Reason in $T$. Suppose
${\cal K}:A\rhd_TB$. Let $J$\/ be the cut of the $T+A$-numbers, which is
isomorphic to
an initial segment of the $T+B$-numbers viewed via $\cal K$.
As above, we may assume $J$\/ to be a $T$-cut.
We have $\Box_T(\Box_TC\to\Box_T\Box_T^JC)$ and, hence,
$\Box_T((\Diamond_TA\et\Box_TC)\to \Diamond_T(A\et\Box_T^JC))$. It follows
that:
\begin{eqnarray*}
 \Diamond_TA\et\Box_TC & \rhd_T & \Diamond_T(A\et\Box_T^JC) \\
                         & \rhd_T & A\et\Box_T^JA \\
                         & \rhd_T & B\et\Box_TC
\end{eqnarray*}
The last step, is witnessed by $\cal K$, using the fact that $\cal K$ \/ is
an end-extension
of $J$\/ and the upwards persistence of $\Sigma^0_1$-sentences.

\bigskip\noindent ${\sf M}_0$ can be viewed as a kind of `{\sf M}-ified'
version of {\sf J5}.
First note that we can rewrite {\sf J5} as the equivalent:
\begin{itemize}
\item $\vdash A\rhd B \to \Diamond  A\rhd B$
\end{itemize}
We get ${\sf M}_0$, by plugging the $\Box C$'s into the consequent.

\subsubsection*{The Principle ${\sf W}^\ast$}

${\sf W}^\ast$ is the following principle.
\begin{itemize}
\item[${\sf W}^\ast$] $\vdash A\rhd B \to (B\et\Box C)\rhd (B\et\Box
C\et\Box\neg\,A)$
\end{itemize}
Dick de Jongh has shown that ${\sf W}^\ast$ is equivalent over {\sf IL} with
${\sf W}+{\sf M_0}$.

\subsubsection*{The Principle ${\sf P}_0$}

For some time ${\sf ILWM_0}$ or, if you wish, ${\sf ILW}^\ast$  stood as
the conjectured candidate for
being $\ilvan{{\sf all}}$. Recently, Albert Visser found a new principle
${\sf P_0}$.
\begin{itemize}
\item[${\sf P_0}$] $\vdash A\rhd \Diamond B \to \Box(A\rhd B) $
\end{itemize}

\noindent The discovery of ${\sf P}_0$
will be described in greater detail in subsection~\ref{histoire}.

The arithmetical verification of ${\sf P}_0$ is as follows. Reason in $T$.
Suppose ${\cal K}:A\rhd_T\Diamond B$. Take a suitably large finite
subtheory $T_0$\/ of $T$.
(We can take $T_0$ standardly finite.)
We certainly have ${\cal K}:(T+A)\rhd(T_0+\Diamond_TB)$. Hence,
(a) $\Box_T({\cal K}:(T+A)\rhd(T_0+\Diamond_TB))$.
On the other hand, by interpretation-existence: (b)
$\Box_T((T_0+\Diamond_TB)\rhd(T+B))$. This last argument works, since an
interpretation can be built
in any arithmetical base theory: we could even have taken $T_0={\sf Q}$,
where {\sf Q}
in Robinson's Arithmetic! Finally, composing (a) and (b), we get:
$\Box_T((T+A)\rhd(T+B))$.

\bigskip\noindent
We can view ${\sf P}_0$ as a `{\sf P}-ified' version of {\sf J5}.
First note that we can rewrite {\sf J5} as the equivalent:
\begin{itemize}
\item $\vdash A\rhd \Diamond B \to  A\rhd B$
\end{itemize}
We get ${\sf P}_0$, by putting a box in front of the consequent.

\bigskip\noindent We can now pose a new conjecture:
${\sf ILW^\ast P_0}\stackrel{?}{=}\ilvan{{\sf all}}$.

\section{Modal Semantics}

\subsection{Veltman Semantics}

Frank Veltman, in response to questions of Albert Visser, found a purely
modal Kripke style semantics for interpretability logic. Frank Veltman was
working on
conditionals at that time. However, $\rhd$ is not a proper
conditional. See \cite{velt:logi85}.

Veltman's semantics extends the well-known Kripke semantics
for $\lol$. Recall that an $\lol$-frame is a pair $\tupel{W,R}$ where $W$\/ is
a nonempty set and $R$\/ is a transitive conversely well-founded relation
on $W^2$. An $\lol$-model is a triple $\tupel{W,R,\Vdash}$
where $\tupel{W,R}$ is an $\lol$-frame and $\Vdash$ is a forcing relation
which commutes in the usual way with the connectives
($w\Vdash A\et B\Leftrightarrow w\Vdash A \mbox{ and } w\Vdash B$, etc.)
and, moreover,
$w\Vdash \Box A
\Leftrightarrow \forall v\; (wRv
\Rightarrow v\Vdash A)$. An $\il$-frame or Veltman frame is a triple
$\tupel{W,R,\{ S_w\mid w\in W\} }$ such that:

\begin{enumerate}
\item
$\tupel{W,R}$ is an $\lol$-frame.
\item \label{sub}
$S_w\subseteq w{\uparrow} \times w{\uparrow}$ ($w{\uparrow} := \{x \in W \mid
wRx \}$).
\item \label{Rin}
$(R\restriction( w{\uparrow} ))\subseteq S_w$.
\item \label{refl}
$S_w$ is reflexive.
\item \label{trans}
$S_w$ is transitive.
\end{enumerate}

\noindent
A Veltman model is a quadruple $\tupel{ W,R,\{ S_w\mid w\in W\} ,\Vdash }$.
Here the triple $\tupel{ W,R,\{ S_w\mid w\in W\} }$ is a Veltman frame and
$\Vdash$ is
a forcing relation with the extra condition that
\begin{itemize}\item
$w\Vdash A\rhd B\Leftrightarrow \forall u\;(w\,R\,u\Vdash
A\Rightarrow \exists v\;\;u\,S_w\,v\Vdash B)$.\\
(We write e.g.\  `$u\,S_w\,v\Vdash B$' for `$uS_wv$\/ and $v\Vdash B$'.)
\end{itemize}

\noindent
Veltman semantics is designed so that $\il$ is sound and complete with
respect to it.

\subsection{Frames}

Consider a frame ${\cal F}=\tupel{ W,R,\{ S_w\mid w\in W\} }$.
We define:
\begin{itemize}
\item ${\cal F}\models A :\Leftrightarrow \mbox{for all forcing relations
$\Vdash$,
   and for all $w{\in}W$,}\;
           w\Vdash A$.
\item
$\cal F$\/ is an $\ilw$-frame if, for any $x$, $ R;S_x$\/
 is conversely well-founded. \\
Here $u(R;S_x)v$\/ if, for some $w$, $uRwS_xv$.
\item
$\cal F$\/ is an
 an $\ilmo$-frame if
$xRyRzS_xuRv\Rightarrow yRv$.
\item
$\cal F$\/ is an $\ilw^*$-frame if it is both an
$\ilw$ and an $\ilmo$-frame.
\item
$\cal F$\/ is an $\ilm$-frame if $yS_xzRu\Rightarrow yRu$.
\item
$\cal F$\/ is an $\ilp$-frame if $xRyRzS_xu\Rightarrow zS_{y}u$.
\end{itemize}
 We have the following
correspondences: $\cal F$\/ is an $\ilw$-frame, an $\ilmo$-frame, an
$\ilw^\ast$-frame,
an $\ilm$-frame, an $\ilp$-frame if, respectively, ${\cal F}\models{\sf W}$,
${\cal F}\models{\sf M}_0$, ${\cal F}\models {\sf W}^\ast$, ${\cal
F}\models {\sf M}$,
${\cal F}\models{\sf P}$.

\subsection{Completeness Results}

The logics $\il$, $\ilm$, $\ilp$ are all
modally complete with respect to their corresponding
classes of frames. (See e.g.\ \cite{japa:logi98} or
\cite{viss:over98}.) In \cite{dejo:comp99} it is shown that $\ilw$ is
also modally complete. $\il$, $\ilw$, $\ilm$ and $\ilp$ can be all shown to
have the finite model property. It follows that they are decidable.
 In \cite{joos:inte98} the modal
completeness of $\ilmo$ is proved. Although this theory is conjectured to be
decidable too, its decidability is still open. Also, the question of the
modal completeness for
${\sf ILW}^\ast$ remains open.

The  arithmetical completeness of $\ilm$ and $\ilp$ was proved by embedding
(the algebras associated with) the Veltman models for $\ilm$, respectively
$\ilp$
into the arithmetical theories. Thus the proofs of the arithmetical
completeness theorem
essentialy involved all three features: modal systems, Veltman semantics
and arithmetical
semantics.

\subsection{The Story of ${\sf P}_0$}\label{histoire}

During fall of $1998$, the progress in developing  the modal completeness
proof of $\ilmo$ stagnated.  It was thought that, perhaps,
it would simplify things if we could strengthen the logic.
Albert Visser tried
to strengthen the frame condition of $\ilmo$ to arrive at a stronger
principle. Remember that the frame condition of $\ilmo$ is:
\begin{itemize}
\item
$xRyRzS_{x}uRv \rightarrow yRv$.
\end{itemize}

\noindent
Instead of demanding an $R$-relation between $x$\/ and $v$, one can
demand an $S_y$-connection between $z$\/ and $v$. If we have
$zS_yv$, we must also have $yRv$, so indeed this move results in
strengthening the frame
condition. A corresponding principle, baptized ${\sf P}_0$, turns out to be:
\begin{itemize}
\item[${\sf P}_0$] $\vdash A \rhd \Diamond B \rightarrow \Box ( A \rhd B )$.
\end{itemize}

\noindent
Clearly every $\ilpo$-frame is an $\ilmo$-frame. If the logic
$\ilpo$ were modal\-ly complete then we would have:
$\ilpo \vdash M_0$ (i.e.\ $\ilpo$
proves every instance of ${\sf M}_0$.) In \cite{joos:inte98} it is
shown that $\ilpo \nvdash {\sf M}_0$ and hence that $\ilpo$ is modally
incomplete.
The proof makes essential use
of $\lver$-models which are a refinement of Veltman models,
invented by Dick de Jongh. All logics
are also sound w.r.t.\ the $\lver$-models, but more distinctions between
principles become
apparent. The main idea is that $S$-relations don't run to a single
world but to a set of worlds. More details can be found in
\cite{joos:inte98}.

The real surprise was that the
 principle $P_0$ ---which came from purely modal considerations---
 is valid in any reasonable arithmetical theory and
hence should be in the core logic $\ilg$.

\subsection{New Principles by Modal Refinements}

If we are looking for principles in $\ilg$, we know for sure that they
should be both in $\ilp$ and in $\ilm$. A priori, there is an infinite
search space but, Veltman models provide pretty good
guidance in this quest. We shall make a convention
on visualizing frame conditions. First, we do not represent all
the relations in the pictures. If $aRb$ and $bRc$ are drawn, we will
rather not draw the  $aRc$\/ that is dictated by transivity.
So by a picture we
actually mean its closure w.r.t.\ the closure conditions
for Veltman models. Secondly, the $R$-relations will be drawn as
straight lines and the $S$-relations as curved lines.\\

\[
\psfig{figure=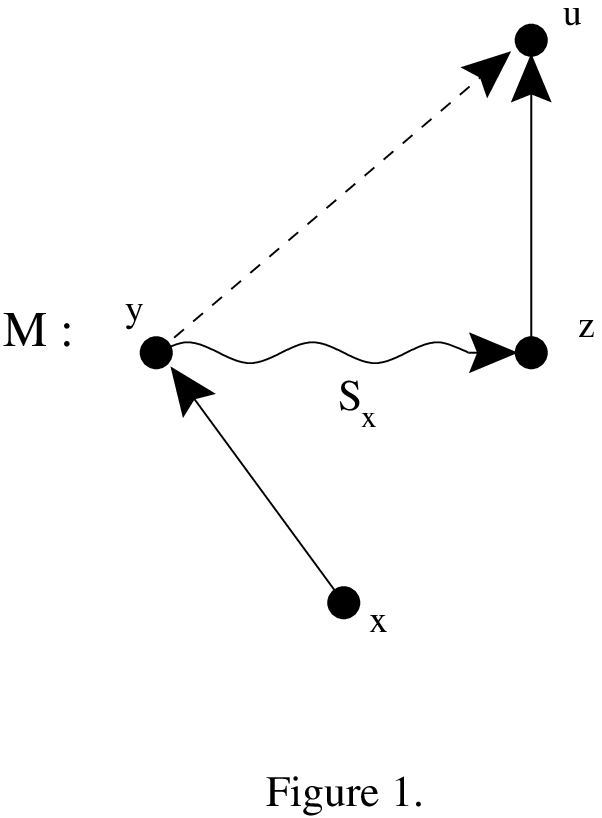,
height=5.5cm}
\]
If some specific relation is imposed by the principle whose frame
condition we want to represent, we will indicate this by drawing a
dashed line.  Bearing this in mind we can visualize the frame condition
of {\sf M}. This condition was:
\[
xRyS_xzRu\Rightarrow yRu.\\
\]

\noindent
In picture 1.\ this condition is represented. The imposed $yRu$\/ is
drawn as a dashed straight line.

\[
\psfig{figure=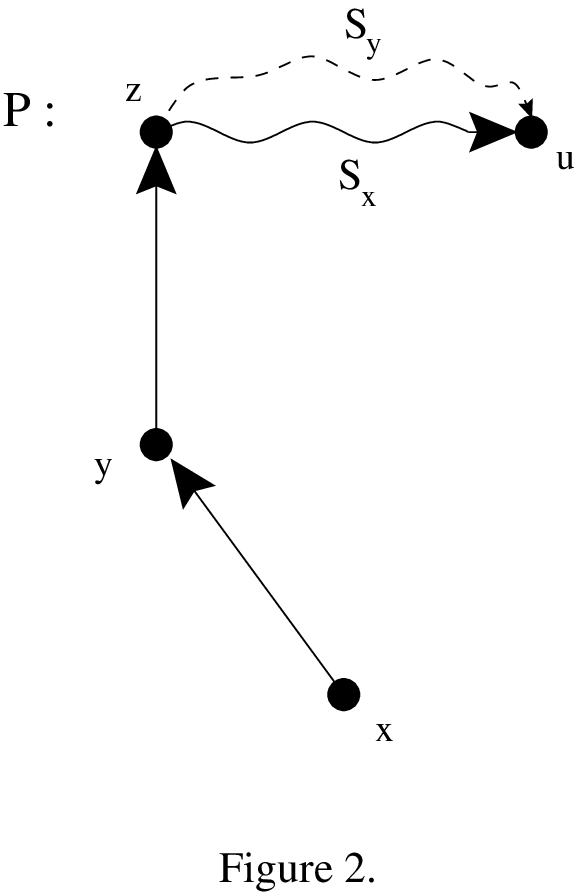
,height=6cm}
\]
When depicting the frame condition of {\sf P}, we get a similar
picture. The frame condition of {\sf P} was:
\[
xRyRzS_xu\Rightarrow zS_{y}u.
\]

\noindent The imposed $zS_{y}u$ is drawn as a curved intersected
arrow. Note that the relations $xRu$\/ and $yRu$\/ are not shown although
they have to exist.

A modal principle in $\ilg$ should hold on all $\ilm$- and
$\ilp$-frames. Consequently the frame condition of this
principle should hold in both frame classes too. This was, of course,
the case for all the
principles considered so far. For example in both
$\ilm$-frames and $\ilp$-frames we have that $R;S_x$ is
conversely well-founded, for any $x$. And this was precisely the frame
condition of
{\sf W}, a principle that holds in any reasonable arithmetical theory.

We
can use the pictorial heuristic to guess new principles. The search space
for new principles is thus confined to principles whose corresponding
frame conditions are shared consequences of both the respective
frame conditions
of {\sf P} and {\sf M}. An example clarifies this concept.\\

\[
\psfig{figure=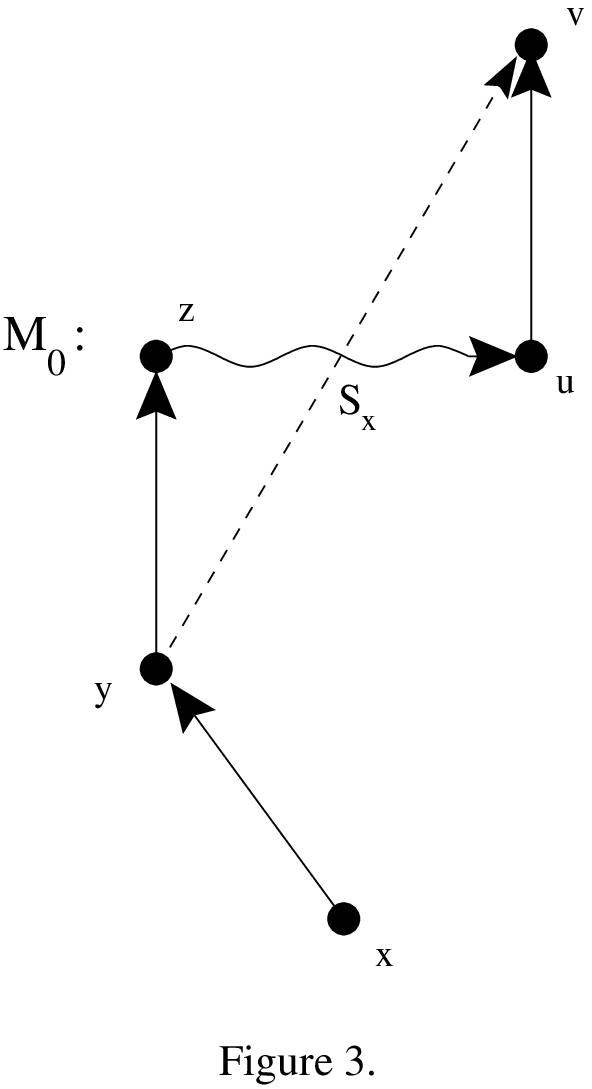,
height=7cm}%
\]
A frame condition is drawn in picture 3. We assume
$xRyRzS_xuRv$. Now the frame condition is that
we impose $yRv$. The corresponding principle is readily found. It is
${\sf M}_0\;\; \vdash A\rhd B\rightarrow
(\Diamond A \wedge \Box C)\rhd (B\wedge \Box C)$.
The $R$-relation between $y$\/ and $v$\/ is implied by the $\ilm$
frame condition because in an $\ilm$ frame we should have $zRv$ and
thus $yRv$ as well.
It is also implied by the frame condition of $\ilp$ because in an
$\ilp$-frame one has $zS_{y}u$ and obviously also $yRu$. And this
again yields $yRv$. The relation $yRv$ is both in the closure of
the frame under $\sf M$ and under $\sf P$. So $yRv$ is in the
intersection.

\[
\psfig{figure=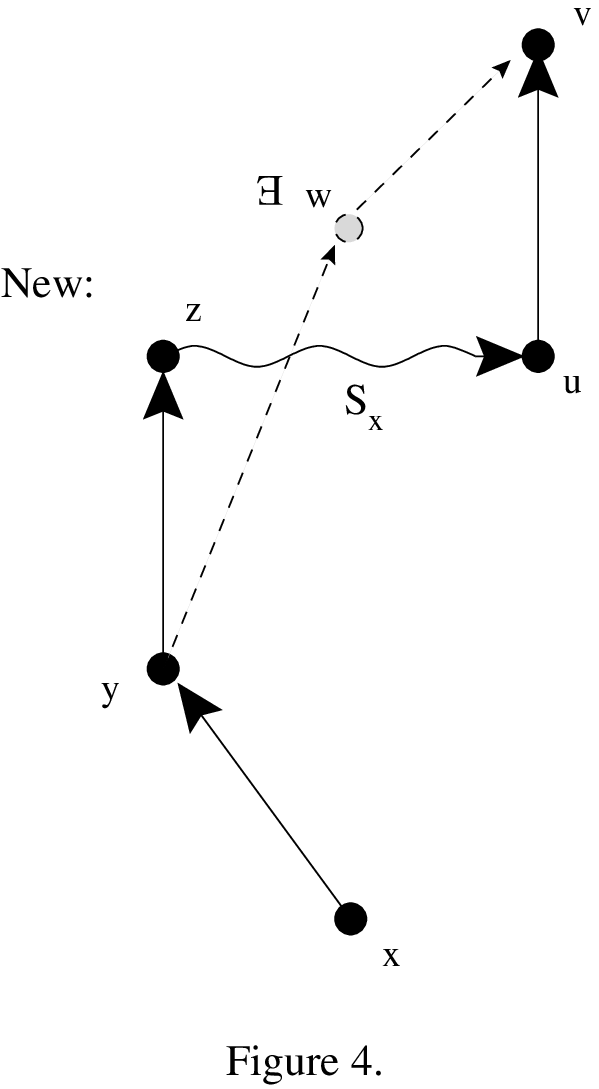
,height=7cm}
\]
Reflection on the previous reasoning tells us that in both the $\ilm$- as in
the $\ilp$-closure it is possible to go from $y$\/ in two $R$-steps to
$v$. In the $\sf M$-case this is $yRzRv$\/ and in the $\sf P$-case this is
$yRuRv$. This idea  could be captured by a somewhat different frame
condition which demands the existence of an intermediate world $w$\/
between $y$\/ and $v$. This condition is represented in figure 4. So
the frame condition is $xRyRzS_xuRv\Rightarrow \exists w \;
yRwRv$ and a corresponding principle is
$\vdash A\rhd B\rightarrow (\Diamond A\wedge \Box \Box C)\rhd (B\wedge \Box
C)$.
In this way we discover a new principle, say ${\sf M}_1$, which is
clearly stronger than ${\sf M}_0$.
{\em It is an open question whether ${\sf M}_1$ is valid in
all reasonable arithmetical theories!}

Another example of a principle in the intersection of
$\ilm$ and $\ilp$ is the principle
$A\rhd \Diamond B \Rightarrow \Box(A\rhd\Diamond B)$.
In \cite{viss:over98} it is shown that this principle is  {\em not}\/
valid in all reasonable arithmetical theories.

\section{Concluding Remarks}


\subsection{The Current Situation}

At the moment of writing, we have a good picture  of the relationships
of the salient logics produced by our quest for $\ilg$.
There are the systems $\ilw$
and $\ilmo$ which are both modally complete and
in the core logic we are looking for. Then there is the
logic $\ilw^*$ which is the union of these two logics. It is
conjectured to be decidable and complete but the problem is still open.
The logic $\ilpo$ is completely independent from $\ilmo$, $\ilw$
and $\ilw^*$. It is also in $\ilg$. The union of all these logics,
 $\ilpow$, is conjectured to be $\ilg$.

\hspace{-0.4cm}
\begin{picture}(0,0)%
\epsfig{file=nieuw.pstex}%
\end{picture}%
\setlength{\unitlength}{2842sp}%
\begingroup\makeatletter\ifx\SetFigFont\undefined%
\gdef\SetFigFont#1#2#3#4#5{%
  \reset@font\fontsize{#1}{#2pt}%
  \fontfamily{#3}\fontseries{#4}\fontshape{#5}%
  \selectfont}%
\fi\endgroup%
\begin{picture}(7575,3090)(901,-2686)
\put(1021,-1111){\makebox(0,0)[lb]{\smash{\SetFigFont{9}{10.8}{\familydefault}{\mddefault}{\updefault}
$\il$
}}}
\put(2446,-1111){\makebox(0,0)[lb]{\smash{\SetFigFont{9}{10.8}{\familydefault}{\mddefault}{\updefault}
$\ilmo$
}}}
\put(3871,-1111){\makebox(0,0)[lb]{\smash{\SetFigFont{9}{10.8}{\familydefault}{\mddefault}{\updefault}
$\ilw^*$
}}}
\put(6346,-1096){\makebox(0,0)[lb]{\smash{\SetFigFont{9}{10.8}{\familydefault}{\mddefault}{\updefault}
$\ilg$
}}}
\put(5101,-1111){\makebox(0,0)[lb]{\smash{\SetFigFont{9}{10.8}{\familydefault}{\mddefault}{\updefault}
$\ilpow$
}}}
\put(8026,-1936){\makebox(0,0)[lb]{\smash{\SetFigFont{9}{10.8}{\familydefault}{\mddefault}{\updefault}
$\ilp$
}}}
\put(7951,-586){\makebox(0,0)[lb]{\smash{\SetFigFont{9}{10.8}{\familydefault}{\mddefault}{\updefault}
$\ilm$
}}}
\put(3751,-1936){\makebox(0,0)[lb]{\smash{\SetFigFont{9}{10.8}{\familydefault}{\mddefault}{\updefault}
$\ilpo$
}}}
\put(2551,-511){\makebox(0,0)[lb]{\smash{\SetFigFont{9}{10.8}{\familydefault}{\mddefault}{\updefault}
$\ilw$
}}}
\put(6001,-811){\makebox(0,0)[lb]{\smash{\SetFigFont{9}{10.8}{\familydefault}{\mddefault}{\updefault}
=
}}}
\end{picture}

\subsection{Two Questions}

We end our paper by formulating two questions of more restricted scope
than our great problem.

\subsubsection*{Problem 1}

The logic $\ilg$ is in the intersection of $\ilm$ and
$\ilp$. But it cannot be equal to this system since e.g.\
 the principle
$A\rhd \Diamond B\rightarrow \Box (A\rhd \Diamond B)$
which is in the intersection,
 is not generally valid. A proof of this  fact is given in
\cite{viss:over98}. The proof employs a heavy result
due to Shavrukov, see
\cite{shav:inter97}. Is there a more direct and more perspicuous proof
of this fact?

\subsubsection*{Problem 2}

Is the principle $A\rhd B\rightarrow (\Diamond A\wedge \Box \Box C)\rhd
(B\wedge \Box C)$\/ arithmetically valid?  It is certainly
in $\ilm$ and $\ilp$ and thus valid both in essentially reflexive and
in finitely axiomatized reasonable arithmetical theories.
Moreover it can be shown to be
valid for $I\Delta_0+B\Sigma_1+\Omega_1$. Yet it is hard to see
why it should be generally valid. In fact we conjecture that it is
not.


\section*{Acknowledgements}
We thank Dick de Jongh for many enlightening conversations. We thank
Rosalie Iemhoff
for her careful reading of the penultimate draft.

\end{document}